\documentclass[]{amsart}
\usepackage{latexsym}
\usepackage{amsmath}
\usepackage{amsfonts}
\usepackage{amsthm}
\usepackage{amssymb}
\usepackage[T1]{fontenc}
\usepackage{hyperref}
\title{A note on two functional inequalities}
\author{W\l{}odzimierz Fechner}
\address{Institute of Mathematics, Lodz University of Technology, ul. W\'olcza\'nska 215, 90-924 \L\'od\'z, Poland}
\email{wlodzimierz.fechner@p.lodz.pl}
\newtheorem{thm}{Theorem}
\theoremstyle{remark}
\theoremstyle{definition}

\newcommand{\R}{\mathbb{R}}

\newcommand{\N}{\mathbb{N}}

\keywords{composite functional inequality}
\subjclass[2010]{Primary: 39B62}

\begin{document}

\begin{abstract}
Motivated by some recent works on the following composite functional equation:
$$f(x^2+yf(z)) =  xf(x)+zf(y), \quad x, y, z \in \R,$$
we solve two related functional inequalities. In each case the solutions are only the constant zero function and the identity on $\R$. 
\end{abstract}

\maketitle 

\section{Introduction}

Composite functional equation 
$$
f(x^2+yf(z)) =  xf(x)+zf(y), \quad x, y, z \in \R\eqno{(E)}\label{E}
$$
appeared in a discussion on Mathematics Stack Exchange \cite{M} (we do not claim that this was its first appearance). Recently, its stability and superstability was studied by A. Najati, M. A. Tareeghee, T. M. Rassias \cite{N} and by D. M. D\u{a}ianu \cite{D}. This research motivates a question on solutions of related inequalities: what will happen if we replace the equality sign with an inequality? An inspection of the reasoning from \cite{M}, \cite{N} and \cite{D} shows that the argumentations used for equations or for a stability problem collapses at some point when applied to an inequality.
 Since the equation is of composite type, and thus nonlinear, then both inequalities need not be fully symmetric. In the present paper we provide the full characterization of both inequalities in the class of real-to-real mappings. It turns out that the sets of solutions are identical and coincide with the solutions of functional equation. No regularity assumptions are needed.

Our main result covers both inequalities, moreover some parts of the proof are analogous.

\begin{thm}
The function $f\colon\R\to\R$ solves at least one of the following inequalities
$$
 f(x^2+yf(z) \geq xf(x)+zf(y), \quad x,y,z\in\R \eqno{(I_+)}\label{P}
$$
$$
  f(x^2+yf(z))\leq xf(x)+zf(y), \quad x,y,z\in\R \eqno{(I_-)}\label{M}
$$
if and only if $f=0$ on $\R$ or $f(x) = x$ for each $x \in \R$.
\end{thm}
\begin{proof}
The "if" part is obvious. We will prove the "only if" part. The proof is divided into a few claims. We assume that $f\colon\R\to\R$ solves $(I_+)$ or $(I_-)$ and is not identically equal to zero. Thus, we aim to show that in both cases $f$ is the identity mapping.

\medskip

\textbf{Claim 1.} \emph{$f(0)=0$, $f(1) = 1$, $f(-1) = -1$ and $f(t) \neq 0$ for $t \neq 0$.} 

\medskip

Substitution $x=y0$ in $(I_+)$ leads to $f(0) \geq zf(0)$ for all $z \in \R$ and to reverse inequality in $(I_-)$. In each case we conclude that $f(0)=0$. Next, substitute $y=z=0$ and $x=0$ to get two alternative pairs of inequalities
\begin{align}
 f(x^2)&\geq xf(x), \quad  f(yf(z))\geq zf(y),  \quad x,y,z\in\R\label{1}\\
 f(x^2)&\leq xf(x), \quad f(yf(z))\leq zf(y), \quad x,y,z\in\R\label{2}.
\end{align}
 \eqref{1} holds if $(I_+)$ is assumed and \eqref{2} is implied by $(I_-)$. 

Let us consider the case $(I_+)$ first. We will show that $f$ attain both a positive and a negative value. 
Assume for a contrary that $f\leq 0$ on $\R$. Then for $x<0$ we get
$$
 0\geq f(x^2)\geq xf(x)\geq 0,
$$
thus $f(x^2)=f(x)=0$, so $f=0$ on $\R$.
If $f\geq 0$ and $f\neq 0$, choose $y_0\in \R$ such that $f(y_0)>0$. Then, for each $z>0$ we have
$
 f(y_0f(z))\geq zf(y_0)>0,
$
which jointly with already established $f(0)=0$ proves that $f(z)>0$, in particular $c:=f(1)>0$. Fix
$x>0$ and put $z=1$ and $y=-x^2/c$ in  $(I_+)$ to get
$0\ge xf(x)+f(-x^2/c)>0$; a contradiction.

Assume now that $f(t)=0$ for some $t\neq 0$. We have $0\ge tf(y)$ for each $y\in \R$; a contradiction with the previous observation on attaining both signs by $f$. In particular
$c=f(1) \neq 0$.
Suppose $c\neq 1$, and substitute $z=1$ and $y=x^2/(1-c)$.  We obtain $xf(x)\leq 0$ for all $x\in \R$. Since $f$ has no zeros except of $f(0)=0$, then $f(x)<0$ for all $x>0$ and consequently $c<0$. 
The first part of \eqref{1} and $f(x^2)<0$ imply $f(x)>0$ for $x<0$. Substitute then $y=-1$ and $z=1$ in the second inequality of \eqref{1} to get $f(-c)\geq f(-1)>0$.
Since $-c>0$, then $f(-c)<0$, which leads to a contradiction. We thus proved that $f(1)=1$. 

Finally apply \eqref{1} for $x=-1$ and the original inequality for $(x,y,z)=(1,-1,1)$ to get
$f(-1)\geq -1$ and $0\geq 1+f(-1)$. Both inequalities together imply that $f(-1)=-1$.

The case $(I_-)$ is a bit different. Put $c=f(1)$ and suppose that $c=0$. Substitute  $z=1$ in \eqref{2} to get $f(y)\geq 0$ for all $y\in \R$. From this we see that for each $x<0$ we get
$
 0\leq f(x^2)\leq xf(x)\leq 0,
$
so $f$ must be constantly equal to zero, which is a contradiction with our assumption. Thus $c\neq 0$.
Suppose then that $c\neq 1$. Substitute $z=1$ and $y=x^2/(1-c)$ in $(I_-)$  to get $xf(x)\geq 0$ for all $x\in \R$, so $c>0$.

Put $z=1$ into the second inequality of \eqref{2}:
\begin{equation}\label{ct}
 f(ct)\leq f(t), \quad t\in\R.
\end{equation}

If $c>1$, then take $b:=1/c\in(0,1)$ and use \eqref{ct} repeatedly to arrive at $f(b^n)\geq c$ for all $n\in \N$.
On the other hand, the firsc inequality of \eqref{2} applied for $b,b^2,b^4,\ldots$
gives$$
 f(b^{2^k})\le b^{2^k-1}f(b), \quad k \in \N.
$$
The right-hand side tends to zero, which leads to a contradiction.

If $0 < c < 1$, then \eqref{ct} gives $f(c^n)\leq c$ and $f(-c^m)\leq f(-1)$ for all $n, m \in \N$.
Put $x=-1$ in \eqref{2} to get $f(-1)\leq -c<0$. Next, substitute in the main inequality
 $x=c^n$, $y=-c^{2n-1}$ and $z=1$ to get
$$
 0\leq c^nf(c^n)+f(-c^{2n-1}) \leq c^{n+1}+f(-1)<0
$$
for $n\in \N$ large enough, which is again a contradiction. We thus proved that $c=1$. 

Verification that $f(-1)=-1$ is similar as before: apply \eqref{2} for $x=-1$ and $(I_-)$ for $(x,y,z)=(1,-1,1)$ to get
$f(-1)\leq -1$ and $f(-1)\geq -1$, so $f(-1)=-1$.

 If $f(t)=0$ for some $t\neq 0$, then $0\le tf(y)$ for all $y\in \R$, which is impossible with $y=1$ or with $y=-1$.

\medskip

\textbf{Claim 2.} \emph{ 
Assume $(I_+)$. Then for each $z> 0$ one has $|f(z)|\leq z$. \\ Assume $(I_-)$. Then for each $z> 0$ one has $|f(z)|\geq z$ and  $f(f(z)^2)>0$. }

\medskip

Proof of this claim is similar in both cases. Assume that  $(I_-)$ holds. Then substitution $z=1$  together with Claim 1 imply
\begin{equation}
f(x^2+y)\leq xf(x) + f(y), \quad x, y \in \R.
\label{4}
\end{equation}
Now fix $z>0$ and denote $a:=f(z)\neq 0$. Use \eqref{2} iteratively to get
$$f(ya^m) \leq z^m f(y), \quad m\in\N,\ y\in\R.$$
Take $m=2$ and $y=1$ and then $m=4$ and $y=-1$:
$$
 f(a^2) \leq  z^2,\qquad f(-a^4) \leq -z^4.
$$
Use \eqref{4} for $x=a^2>0$ and
$y=-a^4$ to get
$$
 0=f(0)\leq a^2f(a^2)+f(-a^4)\leq a^2z^2-z^4.
$$
Since $z>0$, then
$a^2\geq z^2$. From the last estimate and from $f(-a^4)\le-z^4$
we get
$$
 a^2f(a^2)\geq -f(-a^4)\ge z^4>0,
$$
thus $f(f(z)^2)=f(a^2)>0$. 

\medskip

\textbf{Claim 3.} \emph{For each $z >0$ from $f(z)>0$ it follows that $f(z) = z$. } 

\medskip

Denote $a=f(z)>0$. The case of the first inequality is easier. Use \eqref{1} with $y=1$ to get
$f(a) \geq z$. From Claim 2 we have $f(a)\leq a$. Use Claim 2 again to get the reverse inequality $z\leq a$.

Assume that $(I_-)$ holds. Use \eqref{2} with $y=1$ and with $y=-1$ twice to arrive at
$f(a)\leq z$ and $f(-a^2)\le-z^2$.
Apply \eqref{4} with $x=a$ and $y=-a^2$:
$$
 0\le af(a)+f(-a^2)\le af(a)-z^2.
$$
Therefore $f(a)\ge z^2/a>0$. By Claim 2 $f(a)\geq a$, so $a\leq z$.
Use again Claim 2 for $z$ and $a=f(z)>0$ to obtain the reverse inequality.

\medskip

\textbf{Claim 4.} \emph{For each $z >0$ we have $f(z)\in\{-z,z\}$.  } 

\medskip

If $f(z)>0$, then use the previous claim. $f(z)=0$ is impossible, since $z\neq 0$. So, assume that $a:=f(z)<0$.
If $(I_+)$ is assumed, then 
$f(a^2)\geq z^2>0$ by \eqref{1}. Claim 2 applied for $a^2>0$ implies
$z^2\leq f(a^2)\leq a^2$ and applied for $z$ gives $a^2\le z^2$.
Consequently $a^2=z^2$, so $a=-z$.

Assume $(I_-)$. Use the last part of Claim 2 to get
$f(a^2)>0$. Now, we can use Claim 3 for $a^2>0$ to get $f(a^2)=a^2$.
From \eqref{2} we see that $a^2=f(a^2)\leq z^2$, and from Claim 2 we get $a^2\ge z^2$.

\medskip

\textbf{Claim 5.} \emph{If $f$ solves $(I_+)$, then $f(x) = x$ for all $x \in \R$.  } 

\medskip

From \eqref{1} applied for $x=1$ we obtain inductively

$$
 f(t+k)\geq f(t)+k, \quad t\in\R, k\in\N.
$$
For $t=0$ we thus have $f(n)\geq n$.
By Claim 2 $f(n)\le n$, so consequently
$f(n)=n$ for each $n \in \N$. 

Next, we exclude the possibility that $f$ attains a negative value at a positive point. In view of Claim 4, this will force the desired form of $f$ for positive arguments. Suppose that for some $z>0$ we have $f(z)<0$.
By the previous claim $f(z)=-z$. apply \eqref{1} for $y=n\in\N$ we get
$$
 f(-nz)=f(nf(z))\geq zf(n)=nz.
$$
Take $k=\lceil nz\rceil$ and $u=k-nz\in[0,1)$.
We have
$$
 f(u)\geq f(-nz)+k\geq nz+k\geq nz.
$$
But Claim 2 implies $f(u)\leq u<1$ whenever $u>0$ (and if $u=0$, then $f(u)=0$).
Choose a positive integer $n>1/z$ to reach a contradiction. Therefore $f(z)=z$ for all $z>0$.

Take $t>0$ and use $(I_+)$ with $z=-1$ and $y=t$ to get $f(-t)\geq -f(t)=-t$.
From the first part of \eqref{1} applied for $x=\sqrt t$, $s=-t$ we deduce
$$
 0=f(0)\geq \sqrt t\,f(\sqrt t)+f(-t)=t+f(-t).
$$
Therefore $f(-t)\leq -t$, which jointly with a previous estimate leads to $f(-t)=-t$.

\medskip

\textbf{Claim 6.} \emph{ If $f$ solves $(I_-)$, then $f(x) = x$ for all $x \in \R$.  } 

\medskip

Our strategy is the same as in the previous claim, but the argumentation is a bit different.
Suppose that for some $z>0$ we have $a:=f(z)<0$. We have $a=-z$ by Claim 4. From the proof of this claim we also have  $f(a^2)=a^2=z^2>0$.
Use \eqref{2} with $x=z$ and with $y=0$ to get
$$
 f(z^2)\le zf(z)=-z^2.
$$
Since $a^2=z^2$, we get a contradiction $f(z^2)=z^2\le-z^2$. We thus showed that $f$ is positive on the positive halfline and Claim 3 gives that $f(z) = z$ for all $z > 0$.

The proof that $f(-t) = -t$ for $t>0$ is similar like in the previous claim.
\end{proof}

\textbf{Declaration of Generative AI and AI-Assisted Technologies in the Writing Process.}
During the preparation of this work, the author used AI for language polishing, grammar correction and improving the exposition. The core research concepts, mathematical formulations and methodology were entirely conceived and executed by the author. After using this tool, the author thoroughly reviewed and edited the content as needed and takes full responsibility for the accuracy, integrity, and originality of the final publication.

\end{document}